\documentclass[letterpaper, 10 pt,conference]{IEEEconf}%
\usepackage{epsfig}
\usepackage{subfigure}
\usepackage{amsmath}
\usepackage{amsfonts}
\usepackage{graphicx}
\usepackage{graphics} 
\usepackage{amssymb}

\newcommand{\myqedhere}{\hfill$\rule{2mm}{2mm}$}

\newcommand{\agar}{{\rm\ if\ }}

\newcommand{\ow}{{\rm\ otherwise }}
\newcommand{\ind}{\mathbf{1}}

\newcommand{\btheta}{\boldsymbol\theta}

\newcounter{asscount}
\newenvironment{assumption}[1][]{\refstepcounter{asscount}\par\noindent   \textbf{Assumption~\theasscount. #1 } \rmfamily  }{\par}
\newcounter{remcount}
\newenvironment{remark}[1][]{\refstepcounter{remcount}\par\noindent   \textit{Remark~\theremcount: #1 } \rmfamily  }{\par}
\newcounter{lemcount}
\newenvironment{lemma}[1][]{\refstepcounter{lemcount}\par   \textbf{Lemma~\thelemcount. #1 } \rmfamily  }{\par}
\newcounter{themcount}
\newenvironment{theorem}[1][]{\refstepcounter{themcount}\par   \textbf{Theorem~\thethemcount. #1 } \rmfamily }{\par}
\newcounter{propcount}

\begin{document}
%\centerfigcaptionstrue \IEEEoverridecommandlockouts
\title{\textbf{Timeout Control in Distributed Systems Using Perturbation Analysis: Multiple Communication Links}}
\author{\textbf{Ali Kebarighotbi} and \textbf{Christos G. Cassandras }%
\thanks{{\footnotesize The authors' work is supported in part by NSF under
Grant EFRI-0735974, by AFOSR under grant FA9550-09-1-0095 and FA9550-12-1-0113, by DOE
under grant DE-FG52-06NA27490, by ONR under grant
N00014-09-1-1051, and by ARO under grant W911NF-11-1-0227.}}\\Division of Systems Engineering\\and Center for Information and Systems Engineering\\[0pt] Boston University\\[0pt] Brookline, MA 02446\\[0pt] \texttt{alik@bu.edu,cgc@bu.edu}}
\maketitle

\begin{abstract}
Timeout control is a simple mechanism used when direct feedback is
either impossible, unreliable, or too costly, as is often the case
in distributed systems. Its effectiveness is determined by a
timeout threshold parameter and our goal is to quantify the effect
of this parameter on the system behavior.   In this paper, we
extend previous results to the case where there are $N$
transmitting nodes making use of a common communication link
bandwidth. After deriving the stochastic hybrid  model for this
problem, we apply Infinitesimal Perturbation Analysis  to find the
derivative estimates of aggregate average  goodput of the system.
We also derive the derivative estimate of the goodput of a
transmitter with respect to its own timeout threshold which can be
used for local and hence, distributed optimization.
\end{abstract}

\section{Introduction}
Timeout control is a simple mechanism used in many systems where
direct feedback is either impossible, unreliable, or too costly.
This is often the case in distributed systems, where remote
components cannot be observed by a controller (or several
distributed controllers) and information is provided over a
network with unreliable links and delays. A timeout event is
scheduled using a \emph{timer} which expires after some
\emph{timeout threshold} parameter. This defines an expected time
by which some other event should occur or a \textquotedblleft
grace period\textquotedblright\ over which some response about the
system state is required. If no information arrives within this
period, a \textquotedblleft timeout event\textquotedblright\
occurs and incurs certain reactions which are an integral part of
the controller. In distributed systems, where usually control
decisions must be made with limited information from remote
components, timeouts provide a key mechanism through which a
controller can infer valuable information about the unobservable
system states. In fact, as pointed out in \cite{Zhang86}, timeouts
are indispensable tools in building up reliable distributed
systems.

This simple reactive control policy based on timeouts has been
used for stabilizing systems ranging from manufacturing to
communication systems \cite{Jain86}, Dynamic Power Management
(DPM) \cite{CaiLu05},\cite{IraShukGup03},\cite{Rong06} and
software systems \cite{CatPit95},\cite{CoolMobSri99} among others.
Despite its wide usage, quantifications of its effect on system
behavior have not yet received the attention they deserve. In
fact, timeout controllers are usually designed based on heuristics
which may lead to poor results; an important example can be found
in communication protocols, especially TCP (see
\cite{AllPaxStev99},\cite{Jain86},\cite{PsTs07},\cite{Zhang86} and
references therein). There is limited work on finding optimal
timeout thresholds. For instance, in \cite{PerAnuNof03}, a single
queueing model is used for an Automatic Guided Vehicle (AGV)
system and shared tester equipment. In DPM where the aim is to
minimize the average power consumption, \cite{IraShukGup03} and
\cite{CaiLu05} propose a timeout control scheme with the aid of
the theory of competitive analysis. Also, in \cite{CaoChen97}, a
Markov process model is used and Infinitesimal Perturbation
Analysis (IPA) \cite{CassLaf06},\cite{Glasserman91} is applied to
calculate the optimal timeout threshold values. In communication
systems, \cite{GlaMonPen07},\cite{KessMan05},\cite{Libman02} have
attempted to find the optimal TCP retransmission timeouts by
making assumptions on the probability of the transmission failure
in the system. Finally, in \cite{HeGo00} optimal web session
timeouts are calculated so as to reduce the probability of falsely
ending a web session in time sensitive web pages. All such
approaches are limited by their reliance on the distributional
information about
the stochastic processes involved.%\begin{figure}[th] \centering
%\includegraphics[scale =.25]{GeneralModel.jpg}\caption{{\footnotesize {Timeout
%controlled distributed system}}}%
%\label{figs:GeneralModel}%
%\end{figure}
%
%Figure \ref{figs:GeneralModel} shows a high-level model for
%timeout-controlled distributed systems. The control is shown as
%the input $u(S,z)$ where $S$ is the information available to the
%controller, consisting of its own state $X_{c}$, observable system
%state $X_{1}$, and their associated histories $\tilde{X}_{c}$ and
%$\tilde{X}_{1}$. The feedback signal $z$ carries information about
%the state of the remote system and is subject to random
%communication delays. The controller adjusts $u(S,z)$ as a result
%of its input $z$ and timeout events that it generates. We view the
%overall controller-system model as a general \emph{Stochastic
%Hybrid System} (SHS), where the \textquotedblleft system\textquotedblright%
%\ may be time-driven, event-driven or hybrid in itself.

In this paper, we consider the timeout control in the context of
Discrete Event Systems (DES), so that the controller output
includes a response to either a timeout event or an event carrying
information about network congestion. However, since stochastic
DES models can be very complicated to analyze, we rely on recent
advances which abstract a DES into a Stochastic Hybrid System
(SHS) and, in particular, the class of \emph{Stochastic Flow
Models} (SFMs). A SFM treats the event rates as stochastic
processes of arbitrary generality except for mild technical
conditions. The emphasis in using SFMs is not in deriving
approximations of performance measures of the underlying DES, but
rather studying sample paths from which one can derive structural
properties which are robust with respect to the abstraction made.
This is the case, for instance, with many performance gradient
estimates which can be obtained through IPA techniques for general
SHS \cite{CassWarPanYao09},\cite{WarAdMel10}. In addition, a
fundamental property of IPA in SFMs (as in DES) is that the
derivative estimates obtained are independent of the probability
laws of the stochastic rate processes and require minimal
information from the observed sample path. This approach has
proved useful in optimizing various performance metrics in serial
networks \cite{SunCassPan04_2}, systems with feedback control
mechanisms \cite{YuCass06}, scheduling problems
\cite{AliCass09},\cite{AliCass11}, and some multi-class models
\cite{CassSunPanWar03},\cite{SunCassPan04_1}. %Our
%goal here is to use a similar approach for online optimal timeout
%control in various systems within the framework of Fig.
%\ref{figs:GeneralModel}.

In \cite{AliCassCDC11}, we set forth a line of research aimed at
quantifying how timeout threshold parameters affect the system
state and ultimately its behavior and performance. We adopted a
communication system model consisting of one transmitter
submitting packets to a network with stochastic processing times.
We showed how a SFM of a timeout-controlled distributed system can
be obtained and then used to optimize the timeout threshold. This
paper extends the results in \cite{AliCassCDC11} to multiple
communication links and hence, extends this problem to broader
multi-class problems with communication delays. Like
\cite{AliCassCDC11}, we directly control the timeout parameters
for network communication performance. Additionally, moving from a
singe class problem in \cite{AliCassCDC11} to multiple classes, we
show how by carefully selecting the processing rate of each class,
the shared communication channel operates according to a First
Come First Serve (FCFS) policy. We also discuss how we can easily
extend the analysis to non-FCFS frameworks. Finally, we also
extend previous work on timeout control by removing any dependence
on distributional information.

This paper is organized as follows: In Section
\ref{sec:timeoutnd.SFM}, we define the SFM of the system. In
Section \ref{sec:performances}, we apply IPA techniques to average
goodput at each transmitter, as well as total goodput of the
system aggregating the rates for all the transmitting nodes. We
conclude with Section \ref{sec:timeoutnd.conclusion}.

\section{Stochastic Flow Model (SFM)} \label{sec:timeoutnd.SFM}
We are interested in systems where $N$ transmitting nodes send
data to their destinations through a shared channel. Associated
with each node $n=1\ldots,N$, is a ``class'' $n$ of data packets
or tasks which are to be processed in the channel and reach their
destinations in a timely manner. Like the case of a single
communication link, after each transmission by node $n$, a timely
acknowledgement (ACK) is expected from the receiving end. For
transmitter node $n$ to operate normally, the \emph{Round-Trip
Transportation} (RTT) time  calculated for each node $n$ should be
less than a \emph{timeout threshold} parameter $\theta_n\geq 0$.
If no ACK is received within $\theta_n$ units of time from the
transmission, a \emph{timeout} occurs causing node $n$ to
retransmit the timed out data (previously sent at $t-\theta_n$) to
the same destination. The control parameter vector is defined as
$\btheta=[\theta_1,\theta_2,\ldots,\theta_N]^{\textsf{T}}$.
\begin{figure}[th] \centering
\includegraphics[scale =.26]{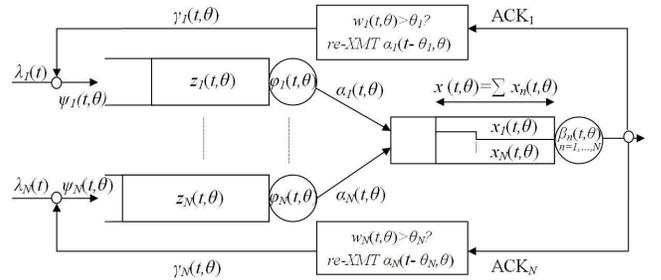}\caption{\footnotesize{The SFM for the timeout system}}%
\label{figs:SFM}%
\end{figure}
We consider a SFM for this problem where $N$ upstream nodes are
competing to get their content through a shared network channel.
The channel latency represents the RTT delay in the underlying
DES.  The SFM for this system is shown in Fig. \ref{figs:SFM} and
is observed over a finite period of time $[0,T]$. Associated with
this system are several nonnegative stochastic processes all
defined on a common probability space $(\Omega,\mathcal{F},P)$.
The transmitters are shown as upstream nodes  receiving exogenous
inflow processes $\{\lambda_n(t)\}$, $n=1,\ldots,N$ and
copy/timeout flow processes $\{\gamma_n(t,\btheta)\}$ intended for
retransmission. The process $\{z_n(t,\btheta)\}$ models the amount
of pending data to be transmitted by node $n$ and evolves
according to
\begin{align}
&\dot{z}_n(t,\btheta)= \frac{dz_n(t,\btheta)}{dt^+}=\nonumber\\
&\left\{\hspace{-2mm}\begin{array}{ll}
0 \hspace{-2mm}& \agar z_n(t,\btheta)=0,\psi_n(t,\btheta)\leq\phi_n(t,\btheta)  \\
\psi_n(t,\btheta)-\phi_n(t,\btheta) \hspace{-2mm}& \ow.
\end{array}\right. \label{zn}
\end{align}
where \begin{equation} \label{psi}
\psi_n(t,\btheta)=\lambda_n(t)+\gamma_n(t,\btheta)\end{equation}
and $\{\phi_n(t,\btheta)\}$ is the \emph{maximal transmission
rate} process at node $n$. Thus, the \emph{actual transmission
rate} from node $n$ to the network buffer is defined as
\begin{align} \alpha_n(t,\btheta)&=\left\{\hspace{-2mm}\begin{array}{ll}
\psi_n(t,\btheta) & \agar z_n(t,\btheta)=0,\psi_n(t,\btheta)\leq\phi_n(t,\btheta) \\
\phi_n(t,\btheta) & \ow \end{array}\right. \label{an}\end{align}

In this paper, we assume that the transmitters have an
\emph{infinite supply} property as follows:
\begin{assumption} \label{assump:timeoutnd.infsupply} $z_n(t,\theta)>0$ for all $t\in[0,T]$. \end{assumption}
By (\ref{an}), this implies
$\alpha_n(t,\btheta)=\phi_n(t,\btheta)>0$ for all $t\in[0,T]$ , so
$\phi_n(t,\theta)$ can henceforth be replaced by
$\alpha_n(t,\theta)$ for the rest of the discussion.

We model the shared network channel with a buffer where the
transmitted fluid from node $n$ accumulates and takes a share
$x_n(t,\btheta)\geq 0$ of the total buffer content
$x(t,\btheta)=\sum_n x_n(t,\btheta)$. The processing rate of the
fluid is governed by nonnegative service processes
$\{\beta_n(t,\btheta)\}$, $n=1,\ldots,N.$

Furthermore, for each $n$, we define the independent process
$\{B_n(t)\}$ as the maximal processing rate of class $n$ when no
other class is being processed by the network resource. Since,
$\beta_n(t,\btheta)$ is the processing rate in the presence of
other fluid classes competing for a share of processing power in
the network, we naturally have $\beta_n(t,\btheta)\leq B_n(t)$.

At any time $t\in[0,T]$, $x_n(t,\btheta)$ in the network follows
the dynamics
\begin{align} &\dot{x}_n(t,\btheta) = \frac{dx_n(t,\btheta)}{dt^+}=\nonumber\\
&\left\{\hspace{-2mm}\begin{array}{ll} 0 \hspace{-3mm}&
\agar x(t,\btheta)=0, \alpha(t,\btheta)\leq\beta(t,\btheta)\\
\alpha_n(t,\btheta) - \beta_n(t,\btheta) \hspace{-3mm}& \ow
\end{array}\right. \label{xndot}
\end{align}
where $\alpha(t,\btheta)=\sum_n\alpha_n(t,\btheta)$ and
$\beta(t,\btheta)=\sum_n\beta_n(t,\btheta)$. We also define the
process
\begin{equation} \label{w_n}
w_n(t,\btheta)=\min_{w\geq 0}\left\{
\int_{t-w}^{t}\beta_n(\tau,\btheta)d\tau =
x_n(t-w,\btheta)\right\}, \end{equation} as the \emph{waiting
time}  process of class $n\in\{1,\ldots,N\}$. We assume that the
value of the processes $\beta_n(t,\btheta)$ and $x_n(t-w,\btheta)$
are known over $t\in[-w,0]$ as the initial pieces of information
required to calculate $w(t,\btheta)$ at $t=0$. Definition
(\ref{w_n}) has a very close relation to the ``time-to-empty''
process as defined in \cite{YaoCass2010},\cite{ChenHuFu2010} since
one way to interpret (\ref{w_n}) is as the earliest time it takes
the server to process the content at $t-w$.

Using the waiting time process, we can define the timed-out
(hence, worthless) portion of class-$n$ fluid in the network
buffer for which $w_n(t,\btheta)>\theta_n$. This portion needs to
be retransmitted by the upstream node $n$. Under normal
conditions, the network generates an acknowledgement flow
indicating a successful transmission previously made at
$t-w_n(t,\btheta)$. In case of a timeout at node $n$ - i.e. the
moment $w_n(t,\btheta)$ exceeds $\theta_n$ - a copy of the fluid
sent at time $t-\theta_n$ is generated and will be put back in
node $n$ for retransmission. This is shown as the feedback flow
process $\gamma_n(t,\btheta)\geq 0$ in Fig. \ref{figs:SFM}. To
summarize, for all $t\in[0,T]$, we have
\begin{equation} \gamma_n(t,\btheta)=\left\{\begin{array}{ll}0 &
\agar w_n(t,\btheta)\leq \theta_n\\
\alpha_n(t-\theta_n,\btheta) & \ow \end{array}\right.
\label{gamman}
\end{equation}
The upstream node treats this copy flow with higher priority and
tries to retransmit it as soon as possible. This is inline with
many retransmission-based communication protocols. Looking at
Fig.\ref{figs:SFM}, we should explicitly show the flow
$\gamma_n(t,\btheta)$ creating contents ahead of those created by
$\lambda_n(t)$. However, we refrained from this to avoid
complicating the illustration.

We look at \emph{average goodput} as the communication performance
metric defined as
\begin{align} J(\btheta,T) &= \mathrm{E}[G(\btheta,T)]\label{J}\\
G(\btheta,T) &= \sum_n\int_0^T
\bigg[\alpha_n(t,\btheta)-2\gamma_n(t,\btheta)\bigg]dt
\label{goodput}
\end{align}
where,  like the one-dimensional problem \cite{AliCassCDC11}, we
use $2\gamma_n(t,\btheta)$ to not only penalize the
retransmissions but also to account for the worthless timed out
fluid in the network buffer which should be processed by the
network resource anyway.

Before proceeding, we make the following assumption on the system
which we also made in \cite{AliCassCDC11}:
\begin{assumption} \label{assump:timeoutnd.noloss} The network buffer is lossless. \end{assumption}
This assumption is merely for simplifying the exposition and can
be removed at the expense of having new state variables tracking
the  loss volume of each class.

\subsection{FCFS Implementation}
The definition of the waiting time (\ref{w_n}) is contingent upon
having a FCFS policy implemented in the network buffer. One of the
differences between a single transmitter model in
\cite{AliCassCDC11} is in how to ensure this policy is preserved.
We achieve this by carefully selecting the network processes.

In general, depending on the processes $\beta_n(t,\btheta)$, if
two fluid particles of class $n$ and $m$, $n\neq m$ are
transmitted at times $\tau_1$ and $\tau_2$, $\tau_2>\tau_1$, if
$\beta_m(t,\btheta)$ is sufficiently larger than
$\beta_n(t,\btheta)$, their order of leaving the network may be
reversed. We can implement both FCFS and non-FCFS policies by
carefully defining the processing rates $\beta_n(t,\btheta)$,
$n=1,\ldots,N.$  We limit the details on the non-FCFS policies to
Remark \ref{rem:nonFCFS} where it is shown that the extension of
the results to the non-FCFS policies is straightforward.

Recall that $\beta_n(t,\btheta)\leq B_n(t)$ since some part of the
resource is consumed by other classes. Moreover, just as in the
DES where a slow service for one type of customer increases the
waiting time of the other customers behind it, if $B_n(t)$ is
small, it not only means slow processing for class $n$, but also a
decreasing effect on other processing rates $\beta_m(t,\btheta)$,
$m\neq n$. The \emph{availability rate} of fluid class $n$ at the
network server can be defined as
\begin{equation}\label{tildalpha} \tilde{\alpha}_n(t,\btheta)=\alpha_n(t-w_n(t,\btheta),\btheta)\quad \forall t\in[0,T]. \end{equation}
%\begin{equation}\label{tildalpha} \tilde{\alpha}_n(t,\btheta)=\alpha_n(t-w_n(t,\btheta),\btheta)
%& \agar n\in I(t)\\
%0 & \ow \end{array}\right. \end{equation}
Accordingly, at any time $t\in[0,T]$, the utilization of the
service associated with class $n=1,\ldots,N$ is as follows:
\begin{equation} \label{rho_n} \rho_n(t,\btheta)
=\frac{\tilde{\alpha}_n(t,\btheta)}{B_n(t)}.
\end{equation}
Hence, at each time $t$, $\rho_n(t,\btheta)$ is a measure of how
engaged the server is with the fluid of class $n$. The following
theorem ensures all these characteristics are delivered by
carefully defining $\beta_n(t,\btheta),\ n=1,\ldots,N$. Compared
to \cite{YaoCass2010},\cite{ChenHuFu2010}, it allows for different
maximal service processes for each class.
\begin{theorem} Let $w\in\mathbb{R}^+$. Assuming the initial conditions $w_n(0,\btheta)\equiv w$ for all
$n$, if
\begin{equation} \beta_n(t,\btheta) =
\frac{\tilde{\alpha}_n(t,\btheta)}{\sum_{m} \rho_m(t,\btheta)} =
\frac{\tilde{\alpha}_n(t,\btheta)}{\sum_{m}
\frac{\tilde{\alpha}_m(t,\btheta)}{B_m(t)}},\quad \forall
t\in[0,T], \label{beta_n}\end{equation} the following statements
are true:
\begin{itemize} \item[$(i)$] There
exists a common waiting time process $w(t,\btheta)$ such that
$w_n(t,\btheta)=w(t,\btheta)$ for all $n$ and all $t\in[0,T]$.
\item[$(ii)$] If $0\leq t_1< t_2$, then $t_1-w(t_1,\btheta)<
t_2-w(t_2,\btheta)$
\end{itemize}\label{lem:FCFS}\end{theorem}
\textit{Proof:} When $x(t,\btheta)=0$, the
buffer clearly operates according to FCFS, so we only consider the
case where $x(t,\btheta)>0$. By proving statement $(i)$, we show
the fluid particles at the head of the queue must have arrived at
the same time; Statement $(ii)$ shows that regardless of the
class, fluid particles leave the system with the
order they have been transmitted.\\
We start by proving statement $(i)$. Differentiating the term in
brackets in (\ref{w_n}) with respect to $t$ reveals
\begin{align}
&\beta_n(t,\btheta)-[1-\dot
w_n(t,\btheta)]\beta_n(t-w_n(t,\btheta),\btheta)=\nonumber\\&[1-\dot
w_n(t,\btheta)]\dot
x_n(t-w_n(t,\btheta),\btheta)\nonumber\end{align} which after
regrouping the terms yields \begin{align} &\dot
w_n(t,\btheta)=1-\frac{\beta_n(t,\btheta)}{\alpha_n(t-w_n(t,\btheta),\btheta)}=1-\frac{\beta_n(t,\btheta)}{\tilde
\alpha_n(t,\btheta)}.\label{wdot}\end{align} Using (\ref{beta_n})
in the last equation gives
\begin{equation}\label{wdotFCFS}\dot{w}_n(t,\btheta)=1-\frac{1}{\sum_m
\rho_m(t,\btheta)},\quad \forall t\in[0,T).\end{equation} for
which we consider the constraint $w_n(t)\geq 0$. Notice that the
RHS is independent of $n$. Hence, with the initial condition that
$w_n(0,\btheta)\equiv w$ for all $n$,
$w_n(t,\btheta)=w_m(t,\btheta)$ for any pair $n,m,\ n\neq m$ and
the FCFS operation when $w_n(t,\btheta)$ is defined for all $n$.
Finally, we can define the waiting time dynamics of a fluid
differential at the head of the network buffer as
%If $\tilde \alpha_n(t,\btheta)$ becomes zero for a non-empty
%interval $I=(\tau_1,\tau_2)$, $w_n(t,\btheta)$ is undefined for
%any $t\in I$. However, because the service for the other classes
%follows $\beta_n(t,\btheta)$ in the theorem statement, at
%$t=\tau_2$, $w_n(\tau_2,\btheta)=w_m(\tau_2,\btheta)$ for any
%$m\in M(\tau_2)$. Therefore, starting at $\tau_2$, all classes
%solve (\ref{wdotFCFS}) with the equal initial conditions. and
%$w_n(t,\btheta)=w_m(t,\btheta)$ for any $n,m\in I(t)$.
\begin{equation} \dot{w}(t,\btheta)=\left\{\begin{array}{ll}0 & \agar x(t,\btheta)=0\\
1-\frac{1}{\sum_m\rho_m(t,\btheta)}
&\ow\end{array}\right.\label{wFCFS}\end{equation} This proves part
$(i)$ in the theorem statement.

For part $(ii)$ notice that by definition, a fluid particle which
is at the head of the queue at time $t$ has arrived at
$t-w(t,\btheta)$. Therefore, we only need to show that for any
$t_1,t_2\in[0,T]$ such that $t_1<t_2$, we have
$t_1-w_1(t_1,\btheta)<t_2-w_2(t_2,\btheta)$. We can write
$t_1-w(t_1,\btheta)<t_2-w(t_2,\btheta)$. It follows that
\begin{align}
&t_1-t_2<w(t_1,\btheta)-w(t_2,\btheta)\nonumber \\
&\qquad=w(t_1,\btheta)-\left[w(t_1,\btheta)+\int_{t_1}^{t_2}\dot{w}(\tau,\btheta)d\tau\right]\nonumber\\
&\qquad=-\int_{t_1}^{t_2}\dot{w}(\tau,\btheta)d\tau\nonumber\\
&\qquad= t_1-t_2-\int_{t_1}^{t_2}
-\frac{1}{\sum_m\rho_m(\tau,\btheta)}d\tau\nonumber\end{align}
which, by (\ref{wFCFS}), gives \begin{align}
&0<\int_{t_1}^{t_2}\frac{1}{\sum_m\rho_m(\tau,\btheta)}d\tau.\nonumber\end{align}
Notice that $\frac{1}{\sum_m\rho_m(\tau,\btheta)}>0$. This proves
part $(ii)$ and the theorem.\myqedhere
 With a more complicated proof, we can show the
theorem statements are true even if Assumption
\ref{assump:timeoutnd.infsupply} does not apply. However, we do
not consider it in this paper.

%\begin{remark} If we consider $\sum_m\rho_m(t,\btheta)$ as the total utility
%of the network resource, (\ref{wFCFS}) is inline with the
%intuition because it indicates that $w(t,\btheta)$ increases when
%$\sum_m\rho_m(t,\btheta)>1$ (unstable buffer) and decreases when
%$\sum_m\rho_m(t,\btheta)<1$ (stable buffer), and remains constant,
%otherwise. \end{remark}
\begin{remark}
We can come up with schemes other than FCFS and generalize the SFM
analysis. For example, if by defining
$$\beta_n(t,\btheta)=B_n(t)\frac{\tilde \alpha_n(t,\btheta)}{\sum_{m}\tilde \alpha_m(t,\btheta)},\quad \forall n.$$
we share the resource capacity according to relative availability
of the fluid classes at the server, (\ref{wdot}) becomes
\begin{equation} \label{nonFCFSw}
\dot{w}_n(t)=1-\frac{B_n(t)}{\sum_m\tilde \alpha_m(t,\btheta)},
\quad \forall t\in[0,T].\end{equation} Now, the right-hand side of
this equation is $n$-dependent, which breaks the FCFS rule.
\label{rem:nonFCFS}
\end{remark}

For the rest of this paper, we make the following assumption:
\begin{assumption} $B_n(t)=B(t)$ for all $t\in[0,T]$ and
$n=1,\ldots,N.$ \label{assump:timoutnd.B}\end{assumption} We can
remove this assumption at the expense of more complexity which
diverts the focus from the main purpose of this analysis.
Assumption \ref{assump:timoutnd.B} is not limiting in the present
problem as the communication channels treat the packets from
different sources or classes equally. Moreover, a byproduct of
Assumption \ref{assump:timoutnd.B} is that
\begin{align}\beta(t,\btheta)&=\sum_n \beta_n(t,\btheta)=\sum_n
\frac{\tilde{\alpha}_n(t,\btheta)}{\sum_m
\frac{\tilde{\alpha}_m(t,\btheta)}{B_m(t)}}\nonumber\\&=B(t)\sum_n
\frac{\tilde{\alpha}_n(t,\btheta)}{\sum_m
\tilde{\alpha}_m(t,\btheta)}=B(t),\quad \forall t\in[0,T],
\label{independentbeta}\end{align} which means the total
processing rate of the channel is independent of the timeout rates
chosen.

\subsection{Transmission Control}
We say node $n$ is operating in a \emph{normal period} when
$w(t,\btheta)\leq \theta_n$. In this mode, we assume that the
inflow rates are determined according to the policy $\pi_{1,n}$
which is designed to increase the transmission rate
$\alpha_n(t,\btheta)$. We adopt a second policy $\pi_{2,n}$ which
applies when the node is in the \emph{timeout period}, i.e.,
$w(t,\btheta)>\theta_n$. Policy $\pi_{2,n}$ is reactive and aims
at reducing the transmission rate until the network channel comes
out of the congestion and can satisfies the requirement
$w(t,\btheta)\leq \theta_n$ again. There can be many choices for
the policies $\pi_{1,n}$ and $\pi_{2,n}$. For the sake of
analysis, we choose the policies as follows:
\begin{align} \dot\alpha_n(t,\btheta) &=
f_{\alpha_n}(t,\btheta)=\left\{\begin{array}{ll} r_a & \agar
w(t,\btheta)\leq \theta_n\\
0 & \ow \end{array}\right.,\label{alphadot}\\
\alpha_n(t,\btheta)&\equiv \alpha_{n,\min} \agar
w(t,\btheta)>\theta_n. \nonumber
\end{align} %Figure \ref{figs:alpha_n} shows an example for
%the maximal transmission rate at node $n$.
%\begin{figure}[th] \centering
%\includegraphics[scale =.25]{rate.jpg}\caption{\footnotesize{A sample for $\alpha_n(t,\btheta)$.}}%
%\label{figs:alpha_n}%
%\end{figure}
\subsection{Stochastic Hybrid Model} Viewed as a SHS, we can conceive of the following SFM operation modes
throughout the sample path: We refer to a period over which
$x(t,\btheta)>0$ and $x(t,\btheta)=0$, as a \emph{Non-Empty
Period} (NEP) and \emph{Empty Period} (EP), respectively.
Moreover, we denote the periods over which $w(t,\btheta)>\theta_n$
and $w(t,\btheta)\leq \theta_n$ by TOP$_n$ and NP$_n$,
respectively.

Let $\tau_k$, $k=1,\ldots,K$, be the SFM event times observed in a
sample path of the system over the interval $[0,T]$. We also
define $\tau_0=0$ and $\tau_{K+1}=T$ for notational convenience
and let $e_k$ be the event occurring at $\tau_k$. We are
interested in the following set of events:
\begin{align}
\mathcal{E}=\big\{E_\lambda,E_{B},[x>0],[x=&0],[w>\theta_n],
[w\leq\theta_n],\nonumber\\&[{\tilde
\alpha_n^+\neq\tilde\alpha_n^-}],[\gamma_n^+\neq\gamma_n^-]\big\}.\label{evset}\end{align}
Here, $E_\lambda$ and $E_{B}$ respectively refer to random jump
events in any $\lambda_n(t)$, $n=1,\ldots,N$ and $B(t)$ which are
not affected by $\btheta$. We call the events with this property
\emph{exogenous}. The start and end of a NEP is marked by the
events $[x>0]$ and $[x=0]$, respectively. $[w>\theta_n]$ and
$[w\leq\theta_n]$ are the events defining the start and end of a
TOP$_n$. Since these events are dependent on the system states, we
categorize them as \emph{endogenous} events. When a timeout
$[w>\theta_n]$ event occurs, according to (\ref{alphadot}),
$\alpha_n(t,\btheta)$ drops. This discontinuity will then be
reflected in $\tilde\alpha_n(t,\btheta)$ and consequently, by
\eqref{beta_n}, in $\beta_n(t,\btheta)$. We refer to the event of
$\tilde \alpha_n(t^+,\btheta)\neq \tilde \alpha_n(t^-,\btheta)$ by
$[{\tilde \alpha_n^+\neq\tilde\alpha_n^-}]$. As the last system
event, $[\gamma_n^+\neq\gamma_n^-]$ is an event of discontinuity
in $\gamma_n(t,\btheta)$ strictly inside a TOP. Since the start
and stop of a TOP$_n$ is already defined by $[w>\theta_n]$ and
$[w\leq\theta_n]$, by \eqref{gamman}, $[\gamma_n^+\neq\gamma_n^-]$
only reflects a discontinuity in $\alpha_n(t-\theta_n,\btheta)$
when $w(t,\btheta)>\theta_n$. This directly affects the objective
function (\ref{goodput}) as it is related to $\gamma_n(t,\btheta)$
defined by (\ref{gamman}). We call $[{\tilde
\alpha_n^+\neq\tilde\alpha_n^-}]$ and $[\gamma_n^+\neq\gamma_n^-]$
\emph{induced} events because these events are bound to occur
after the occurrence of a \emph{triggering} event in the past.
Generally, the triggering event can be either an exogenous,
endogenous or itself an induced event. More information on the
induced events can be found in \cite{CassWarPanYao09}. %In our model, any
%$[w>\theta_n]$  at $\tau_{m}$ creates a discontinuity in $\tilde
%\alpha_n(t,\btheta)$ at $\tau_k>\tau_m$, and as such, in the
%resource allocation (see (\ref{beta_n})). Hence, $[{\tilde
%\alpha_n^+\neq\tilde\alpha_n^-}]$ is an induced event which is
%triggered by $[w>\theta_n]$ event. Moreover, $[w>\theta_n]$ also
%causes a jump in $\alpha_n(\tau_k-\theta_n,\btheta)$ at
%$\tau_k=\tau_m+\theta_n$. If we also have
%$w(\tau_k,\btheta)>\theta_n$, by (\ref{gamman}) we get
%$e_k=[\gamma_n^+\neq\gamma_n^-]$.

The delays between $[w>\theta_n]$ and its associated $[{\tilde
\alpha_n^+\neq\tilde\alpha_n^-}]$ and $[\gamma_n^+\neq\gamma_n^-]$
call for new state variables. The role of these state variables is
to provide \emph{timers} which trigger when $[w>\theta_n]$ occurs
and measure the amount of time until the associated induced event.
For each event time $\tau_k$ we define the state variables
\begin{subequations} \label{timer1} \begin{align}
\dot h_{1}(k,t,\btheta)&=\left\{\begin{array}{ll} -\beta(t) &
\agar \exists
n: e_k=[w>\theta_n]\\
0 & \ow \end{array}\right.\forall t\geq \tau_k,\\
h_{1}(k,\tau_k^+,\btheta)&=\left\{\begin{array}{ll} x(t_k,\btheta)
& \agar
e_k=[w>\theta_n]\\
0 & \ow \end{array}\right.
\end{align}\end{subequations}
\begin{subequations} \label{timer2} \begin{align}
\dot h_{2}(k,t,\btheta)&=\left\{\begin{array}{ll} -1 & \agar
\exists
n: e_k=[w>\theta_n]\\
0 & \ow \end{array}\right.\forall t\geq \tau_k,\\
h_{2}(k,\tau_k^+;\btheta)&=\left\{\begin{array}{ll} \theta_n &
\agar
e_k=[w>\theta_n]\\
0 & \ow \end{array}\right.
\end{align}\end{subequations}
with the constraint $h_i(k,t,\btheta)\geq 0$, $i=1,2$,
$k=1,\ldots,K$. To identify the active timers on
$[\tau_k,\tau_{k+1})$, for each $k=1,\ldots,K$, we define the
index set
\begin{equation} \label{triggeringindex} \Phi_{k}=\{m\leq k: \exists
i \mbox{ s.t. } h_{i}(m,t,\btheta)>0\  \forall\
t\in[\tau_k,\tau_{k+1})\}.
\end{equation}

\section{Performance Optimization by IPA}
The following assumption ensures existence of the IPA derivatives:
\begin{assumption} With probability 1, no two events can occur at
the same time unless one causes the other.
\label{assump:timeoutnd.nosimul}\end{assumption}

\label{sec:performances} Considering the objective
(\ref{goodput}), let us define the index set \begin{equation}
\label{omegaset} \Omega_n = \big\{k: w(t,\btheta)>\theta_n\
\forall t\in[\tau_k,\tau_{k+1})\big\}
\end{equation} marking all the event times in TOP$_n$ including its start. The other possible event in $\Omega_n$
is $[\tilde\alpha_n^+\neq\tilde\alpha_n^-]$ which by Assumption
\ref{assump:timeoutnd.nosimul}, cannot occur independently of
$[w>\theta_n]$ and $[w\leq\theta_n]$. Since by (\ref{gamman}),
when $w(t,\btheta)>\theta_n$, we have
$\gamma_n(t,\btheta)=\alpha_n(t-\theta_n,\btheta)$, we can
decompose (\ref{goodput}) as follows:
\begin{subequations}\label{decomposedG}\begin{align} G(T,\btheta)&=\sum_n
G_n(T,\btheta),\label{decomposedG1}\\
G_n(T,\btheta)&=\bigg\{\sum_{k=0}^K\int_{\tau_k}^{\tau_{k+1}}\alpha_n(t,\btheta)dt\nonumber\\&-2\sum_{k\in\Omega_n}\int_{\tau_k}^{\tau_{k+1}}\alpha_n(t-\theta_n,\btheta)dt\bigg\},\
n=1,\ldots,N. \label{decomposedG2} \end{align}\end{subequations}
%To keep the notation manageable, we drop $\btheta$ from function
%arguments except for $J(\btheta)$ and $G(T,\btheta)$.
\subsection{IPA Estimation}
Let us define $\tau_{k,j}^\prime\equiv\frac{\partial
\tau_k}{\partial\theta_j}$ for $k=0,\ldots,K+1$ and
$j=1,\ldots,N$. Also, for a real valued function $f_n(t,\btheta)$
associated with node $n$, we define the partial derivatives
$f^\prime_{n,j}(t,\btheta)\equiv \frac{\partial
f_n(t,\btheta)}{\partial \theta_j}$. Differentiating
$G_n(T,\btheta)$ with respect to $\theta_j$ and noticing by
$\tau_0=0,$ and $\tau_{K+1}=T$, that
$\tau_{0,j}^\prime=\tau_{K+1,j}^\prime=0$ for any $j$ reveals
\begin{align}
&\frac{dG_{n}(T,\btheta)}{d\theta_j}=
\bigg\{\sum_{k=1}^{K}[\alpha_n(\tau_k^-)-\alpha_n(\tau_k^+)]\tau_{k,j}^\prime\nonumber\\
&\ -2\sum_{k\in\Omega_n}\bigg[\tau_{k+1,j}^\prime
\alpha_n([\tau_{k+1}-\theta_n]^-)-\tau_{k,j}^\prime
\alpha_n([\tau_k-\theta_n]^+)\nonumber\\&\
+\sum_{k=0}^{K}\int_{\tau_k}^{\tau_{k+1}}\alpha^\prime_{n,j}(t)dt+\int_{\tau_k}^{\tau_{k+1}}\frac{d
\alpha_n(t-\theta_n)}{d \theta_j}dt\bigg]\bigg\}.\label{Gnprime}
\end{align}
 The expression in the last integral
in (\ref{Gnprime}) can be written as
\begin{align} \frac{d
\alpha_n(t-\theta_n,\btheta)}{d \theta_j}&=\frac{\partial
\alpha_n(\tau,\btheta)}{\partial
\theta_j}\bigg|_{\tau=t-\theta_n}\hspace{-5mm}+\frac{\partial
(t-\theta_n)}{\partial \theta_j}\dot
\alpha_n(t-\theta_n,\btheta)\nonumber\end{align} where
$\frac{\partial \alpha_n(\tau,\btheta)}{\partial
\theta_j}\big|_{\tau=t-\theta_n}\equiv\alpha^\prime_{n,j}(\tau-\theta_n,\btheta)$.
Moreover, according to (\ref{alphadot}),
$\dot{\alpha}_n(t-\theta_n,\btheta)=r_n$ if
$w_n(t-\theta_n,\btheta)<\theta_n$ (i.e., node $n$ not in a
TOP$_n$ at $t-\theta_n$). Noting that $\frac{\partial
(t-\theta_n)}{\partial \theta_j}=-1$ only if $j=n$ and 0,
otherwise, we find \begin{align} \frac{d
\alpha_n(t-\theta_n,\btheta)}{d \theta_j} =&
\alpha^\prime_{n,j}(t-\theta_n,\btheta)\nonumber\\&-\left\{\begin{array}{ll}r_n
&
\agar w(t-\theta_n,\btheta)\leq \theta_n, j=n\\
0 & \ow
\end{array}\right. \label{delayedaprime}
\end{align} Looking at the performance
objective as well as the conditions under which an event is
triggered, the state vector of the system is comprised of
$\alpha_n(t,\btheta)$, $\gamma_n(t,\btheta)$,  $n=1,\ldots, N$,
$x(t,\btheta)$ and $w(t,\btheta)$. However, note that
$\gamma_n(t,\btheta)=\alpha_n(t-\theta_n,\btheta)$ when
$w_n(t,\btheta)>\theta_n$. Thus, its derivative can be obtained
from that of $\alpha_n(t,\btheta)$. Furthermore, Assumption
\ref{assump:timoutnd.B}, allows us to have the following useful
lemma which reduces the number of states to only two:
\begin{lemma} If the network channel operates according to the FCFS policy, we get
$$\frac{\partial w(t,\btheta)}{\partial \theta_j} =  \frac{1}{\tilde\alpha(t,\btheta)}\frac{\partial \tilde{x}(t,\btheta)}{\partial \theta_j},\quad \forall t\in[0,T), j=1,\ldots,N,$$
where $\tilde{x}(t,\btheta)=x(t-w(t,\btheta),\btheta)$ for all
$t\in[0,T]$.\label{lem:wprime}\end{lemma} \textit{Proof:}
Differentiating (\ref{w_n}) with respect to $\theta_j$
$j\in\{1,\ldots,N\}$ gives
\begin{align}
&w^\prime_{n,j}(t,\btheta)\beta_n(t-w_n(t,\btheta),\btheta)+\int_{t-w_n(t,\btheta)}^t
\frac{d\beta_n(\tau,\btheta)}{d\theta_j}d\tau\nonumber\\&=\frac{\partial
x_n(\tau,\btheta)}{\partial
\theta_j}\bigg|_{\tau=t-w_n(t,\btheta)}+\frac{\partial
x_n(\tau,\btheta)}{\partial \tau}\frac{\partial
\tau}{\partial\theta_j}\bigg|_{\tau=t-w_n(t,\btheta)}\nonumber\\&=x^\prime_{n,j}(t-w_n(t,\btheta),\btheta)-\dot
x_n(t-w_n(t,\btheta),\btheta) w_{n,j}^\prime(t,\btheta)\nonumber
\end{align}
where we have $\dot
x_n(t-w_n(t,\btheta),\btheta)=\alpha_n(t-w_n(t,\btheta),\btheta)-\beta_n(t-w_n(t,\btheta),\btheta)$.
Thus,\begin{align}
&w^\prime_{n,j}(t,\btheta)\beta_n(t-w_n(t,\btheta),\btheta)+\int_{t-w_n(t,\btheta)}^t\frac{d\beta_n(\tau,\btheta)}{d\theta_j}d\tau\nonumber\\
&=x^\prime_{n,j}(t-w_n(t,\btheta),\btheta)\nonumber\\&\quad-[\alpha_n(t-w_n(t,\btheta),\btheta)-\beta_n(t-w_n(t,\btheta),\btheta)]w^\prime_{n,j}(t,\btheta)\nonumber
\end{align}
Regrouping terms yields the following for all $t\in[0,T]$:
$$w^\prime_{n,j}(t,\btheta)=\frac{x^\prime_{n,j}(t-w_n(t,\btheta),\btheta)-\int_{t-w_n(t,\btheta)}^t \frac{d\beta_n(\tau,\btheta)}{d\theta_j}d\tau}{\tilde\alpha_n(t,\btheta)}.$$
Using Theorem \ref{lem:FCFS} and by Assumption
\ref{assump:timeoutnd.infsupply}, $w_n(t,\btheta)=w(t,\btheta)$
for all $n=1,\ldots,N$ and all $t\in[0,T)$. Since
$w^\prime_n(t,\btheta)=w^\prime(t,\btheta)$ for all
$n\in\{1,\ldots,N\}$ we can add all the numerators and
denominators and write \begin{align}\frac{\partial
w(t,\btheta)}{\partial \theta_j}=\frac{1}{\sum_n\tilde
\alpha_n(t,\btheta)}\sum_n
\bigg[&x^\prime_{n,j}(t-w(t,\btheta),\btheta)\nonumber\\&-\int_{t-w(t,\btheta)}^t
\frac{d\beta_n(\tau,\btheta)}{d\theta_j}d\tau
\bigg].\nonumber\end{align} Since
$\sum_n\tilde\alpha_n(t,\btheta)=\tilde\alpha(t,\btheta)$ and
$\sum_n
\frac{d\beta_n(t,\btheta)}{d\theta_j}=\frac{d\beta(t)}{d\theta_j}=\frac{dB(t)}{d\theta_j}=0$,
we find that for $j=1,\ldots,N,$
$$\frac{\partial w(t,\btheta)}{\partial \theta_j}= \frac{1}{\tilde\alpha(t,\btheta)}\frac{\partial x(t-w(t,\btheta),\btheta)}{\partial \theta_j}=\frac{1}{\tilde\alpha(t,\btheta)}\frac{\partial \tilde{x}(t,\btheta)}{\partial \theta_j}.$$
\myqedhere
\subsection{IPA equations}
Before proceeding, we provide a brief review of the IPA framework
for general stochastic hybrid systems as presented in
\cite{CassWarPanYao09}. If $s(t,\btheta)\in\mathbb{R}^M$ is the
state vector of the SFM, IPA specifies how changes in $\btheta$
influence $s(t,\btheta)$ and the event times $\tau_{k}$ and,
ultimately, how they influence interesting performance metrics
which are generally expressed in terms of these variables. Let us
assume that over an interval $[\tau_{k},\tau_{k+1})$, the SFM is
at some mode during which the time-driven state satisfies $\dot
{s}=f_{k}(s,\btheta,t)$ for some $f_{k}:\mathbb{R}^{M}\times\mathbb{R}%
^{N}\times\lbrack0,T)\rightarrow\mathbb{R}^{M}$. Let
$s^{\prime}(t)\equiv \frac{\partial
s(t)}{\partial\theta}\in\mathbb{R}^{M}\times\mathbb{R}^{N}$ be the
Jacobian matrix for all state derivatives. It is shown in
\cite{CassWarPanYao09} that, for any
$t\in\lbrack\tau_{k},\tau_{k+1})$, $s^{\prime }(t)$ satisfies:
\begin{equation}
\frac{d}{dt}s^{\prime}(t)=\frac{\partial f_{k}(t)}{\partial
s}s^{\prime
}(t)+\frac{\partial f_{k}(t)}{\partial\theta}\label{{ipacontupdate}}%
\end{equation}
for $t\in\lbrack\tau_{k},\tau_{k+1})$ with boundary condition:
\begin{equation}
s^{\prime}(\tau_{k}^{+})\ =\
s^{\prime}(\tau_{k}^{-})+\tau_{k}^{\prime}\left[
f_{k-1}(\tau_{k}^{-})-f_{k}(\tau_{k}^{+})\right]  \label{ipaxcontdiscupdate}%
\end{equation}
for $k=0,\ldots,K$ if $s(t)$ is continuous at $\tau_k$ and
otherwise,
\begin{equation} s^{\prime}(\tau_{k}^{+})\ =\
r^{\prime}(\tau_{k}). \label{ipaxdiscontdiscupdate}%
\end{equation}

An exogenous event at $\tau_{k}$ is not a function of $\btheta$ so
$\tau_{k,j}^{\prime}=0$ for any $j=1,\ldots,N$. In our model these
include $E_\lambda$ and $E_B$. However, for every endogenous
event, $e_{k}$ at $\tau_{k}$ there exists a continuously
differentiable function
$g_{k}:\mathbb{R}^{n}\times\Theta\rightarrow\mathbb{R}$ such that
$\tau _{k}\ =\ \min\{t>t_{k-1}\ :\ g_{k}\left(
s(t;\btheta),\btheta\right)  =0\}$. It is shown in
\cite{CassWarPanYao09} that
\begin{equation}
\tau_{k}^{\prime}=-\left[  \frac{\partial g_{k}}{\partial s}f_{k}(\tau_{k}%
^{-})\right]  ^{-1}\left(  \frac{\partial g_{k}}{\partial\theta}%
+\frac{\partial g_{k}}{\partial s}s^{\prime}(\tau_{k}^{-})\right)
\label{ipaendoprime}%
\end{equation}
if $e_{k}\in\mathcal{E}$ is endogenous and defined as long as
$\frac{\partial g_{k}}{\partial s}f_{k}(\tau_{k}^{-})\neq0$. In
addition to the exogenous and endogenous events, we also have
induced events. Using the results in \cite{CassWarPanYao09}, in
this case, if $ e_k \mbox{ is induced by } e_{m}$, $m<k$,  we can
write
\begin{align}
\label{induprime}
\frac{\partial\tau_k}{\partial\theta_j}=-\left[\dot
h_i(m,\tau_k^-,\btheta)\right]^{-1}\bigg(&\frac{\partial
s(\tau_k^-)}{\partial \theta_j}^\textsf{T}\frac{\partial
h_i}{\partial s}\nonumber\\&+\frac{\partial
h_i(m,\tau_k^-,\btheta)}{\partial\theta_j}\bigg).\end{align} where
$\frac{\partial h_i}{\partial s}$ is the gradient vector of
partial derivatives $h_i$ with respect to state variables and
depending on whether $[\tilde{\alpha}_n^+\neq\tilde{\alpha}_n^-]$
or $[\gamma_n^+\neq\gamma_n^-]$ has occurred, we use $h_1$ or
$h_2$ as defined in (\ref{timer1}) and (\ref{timer2}).
\subsubsection{Event-time derivatives}
We will determine $\frac{\partial\tau_k}{\partial\theta_j}$ for
each event in $\mathcal{E}$ and for each $\theta_j,\
j=1,\ldots,N$. We exclude $E_\lambda$ and $E_B$ as they are
exogenous events with zero event-time derivative. Recalling
$\alpha(t)=\sum_n \alpha_n(t)$,
$\tilde{\alpha}(t)=\sum_n\tilde\alpha_n(t)$ and $\beta(t)=\sum_n
\beta_n(t)$, the following lemma gives the derivatives for
exogenous and induced events:

\begin{lemma} Under policies $\pi_{1,n}$ and $\pi_{2,n}$, $n=1,\ldots,N$, for any $j=1,\ldots,N$, we have\\
{\begin{tabular}{ll} $(i)$ & If $e_k=[x>0]$: $\
\tau_{k,j}^\prime=
-\frac{1}{\dot\alpha(\tau_k^-)-\dot\beta(\tau_k^-)}\frac{\partial\alpha(\tau_k^-,\btheta)}{\partial
\theta_j}$ \\ \vspace{1mm}
 $(ii)$ & If $e_k=[x=0]$:
$\ \tau_{k,j}^\prime =
-\frac{1}{\alpha(\tau_k^-,\btheta)-\beta(\tau_k^-)}\frac{\partial
x(\tau_k^-,\btheta)}{\partial\theta_j}$\\\vspace{1mm}
 $(iii)$ & If
$e_k=[w>\theta_n]$ or $[w\leq\theta_n]$:\\
&$\tau_{k,j}^\prime=\frac{\ind_{n=j}\tilde
\alpha(\tau_k,\btheta)-\frac{\partial
\tilde{x}(\tau_k^-,\btheta)}{\partial\theta_j}}{\tilde
\alpha(\tau_k^-,\btheta)-{\beta(\tau_k^-)}},$ \\
\vspace{1mm} $(iv)$ & If $e_k=[\tilde{\alpha}_n^+\neq\tilde{\alpha}_n^-]$:\\
&$\tau_{k,j}^\prime=\frac{1}{\beta(\tau_k^-)}\left[\frac{\partial
x(\tau_{m}^+,\btheta)}{\partial\theta_j}+\tau_{m,j}^\prime
\alpha(\tau_m,\btheta)\right]$ \\ \vspace{1mm}
  $(v)$ & If $e_k=[\gamma_n^+\neq\gamma_n^-]$:
$\tau_{k,j}^\prime=\tau_{m,j}^\prime+\ind_{n=j}$\\
\end{tabular}}\\
where $\ind_{n=j}$ is an indicator function being 1 when $n=j$ and
0, otherwise and $\tau_m<\tau_k$ is the time of the triggering
event for the induced event at $\tau_k$.
\label{lem:evtimederCh5}\end{lemma} \textit{Proof:} Starting with
part $(i)$, we can invoke (\ref{ipaendoprime}) with
$g_k(t,\btheta)=\alpha(t,\btheta)-\beta(t)$. Noticing
$\frac{\partial g_k}{\partial \alpha}=1$ and $\frac{\partial
g_k}{\partial \beta}=-1$, we find that $\frac{\partial
g_k}{\partial \alpha}\dot\alpha(\tau_k^-,\btheta)+\frac{\partial
g_k}{\partial
\beta}\dot\beta(\tau_k^-)=\dot\alpha(\tau_k^-,\btheta)-\dot\beta(\tau_k^-).$
Moreover, $\frac{\partial g_k(\tau_k^-,\btheta)}{\partial
\theta_j}=\frac{\partial
\alpha(\tau_k^-,\btheta)}{\partial\theta_j}-\frac{\partial
\beta(\tau_k^-)}{\partial\theta_j}$. However, by
(\ref{independentbeta}) we know that $\beta(t)=B(t)$ and is
independent of the control parameters. Thus, we get
$\frac{\partial g_k(\tau_k^-,\btheta)}{\partial
\theta_j}=\frac{\partial\alpha(\tau_k^-,\btheta)}{\partial\theta_j}$.
Inserting the results in the expression for (\ref{ipaendoprime})
verifies part $(i)$. For part $(ii)$, notice that we have an
endogenous event with $g_k(t,\btheta)=x(t,\btheta)$. Noticing that
when $x(t,\btheta)>0$,
$\dot{x}(t,\btheta)=f_x(t,\btheta)=\alpha(t,\btheta)-\beta(t)$.
Hence, $\frac{\partial g_k}{\partial
x}f_x(\tau_k^-,\btheta)=\alpha(\tau_k^-,\btheta)-\beta(\tau_k^-)$.
Also, $\frac{\partial g_k(\tau_k^-,\btheta)}{\partial
\theta_j}=\frac{\partial x(\tau_k^-,\btheta)}{\partial \theta_j}$.
These give
$\frac{\partial\tau_k}{\partial\theta_j}=-\frac{1}{\alpha(\tau_k^-,\btheta)-\beta(\tau_k^-)}\frac{\partial
x(\tau_k^-,\btheta)}{\partial \theta_j}$ which is exactly the
claim of part $(ii)$. For part $(iii)$, we have
$g_k(\tau_k,\btheta)=w(\tau_k,\btheta)-\theta_n$. By
(\ref{wdotFCFS}), and Assumption \ref{assump:timoutnd.B}, we have
\begin{align} \dot{w}(t,\btheta)=f_w(t,\btheta)=1-\frac{1}{\sum_m
\rho_m(t,\btheta)}&=1-\frac{B(t)}{\tilde\alpha(t,\btheta)}\nonumber\\
&=\frac{\tilde\alpha(t,\btheta)-\beta(t)}{\tilde\alpha(t,\btheta)}.\nonumber
\end{align}
Since $\frac{\partial g_k}{\partial w}=1$, we find $\frac{\partial
g_k}{\partial
w}f_w(\tau_k^-,\btheta)=\frac{\tilde\alpha(\tau_k^-,\btheta)-\beta(\tau_k^-)}{\tilde\alpha(\tau_k^-,\btheta)}$.
Finally, $\frac{\partial g_k(\tau_k^-,\btheta)}{\partial
\theta_j}= \frac{\partial w(\tau_k^-,\btheta)}{\partial
\theta_j}-\frac{\partial \theta_n}{\partial \theta_j}$. Clearly,
$\frac{\partial \theta_n}{\partial \theta_j}=\ind_{n=j}$. By Lemma
\ref{lem:wprime}, we get $$\frac{\partial
w(\tau_k^-,\btheta)}{\partial \theta_j}=
\frac{1}{\tilde\alpha(\tau_k^-,\btheta)}\frac{\partial
x([\tau_k-w(\tau_k,\btheta)]^-,\btheta)}{\partial \theta_j}.$$
Thus, we find that
$$\frac{\partial g_k(\tau_k^-,\btheta)}{\partial
\theta_j}=
\frac{1}{\tilde\alpha(\tau_k^-,\btheta)}\left(\frac{\partial
x([\tau_k-w(\tau_k,\btheta)]^-)}{\partial
\theta_j}-\ind_{n=j}\right).$$ Putting these results in
(\ref{ipaendoprime}) gives
\begin{align} \frac{\partial \tau_k}{\partial\theta_j}&= -\frac{1}{\frac{\partial
g_k}{\partial w}f_w(\tau_k^-,\btheta)}\frac{\partial
g_k(\tau_k^-,\btheta)}{\partial \theta_j}\nonumber\\
&=\frac{\ind_{n=j}\tilde \alpha(\tau_k,\btheta)-\frac{\partial
\tilde{x}(\tau_k^-,\btheta)}{\partial\theta_j}}{\tilde
\alpha(\tau_k^-,\btheta)-{\beta(\tau_k^-)}}\nonumber\end{align}
which completes the proof for part $(iii)$. For case $(iv)$, we
have the condition $h_1(\tau_k,\btheta)=0$ which assuming the
triggering time is $\tau_m$ can be written as
$$h_1(m,\tau_k,\btheta)=x(\tau_m,\btheta)-\int_{\tau_m}^{\tau_k}\beta(t)dt =0.$$ Notice
that, $\dot{h}_1(\tau_k^-,\theta)= -\beta(\tau_k^-)$ and
\begin{align}&\frac{\partial h_1(m,\tau_k^-,\btheta)}{\partial
\theta_j}=\frac{\partial x(\tau_m,\btheta)}{\partial
\btheta_j}+\frac{\partial \tau_m}{\partial
\theta_j}\beta(\tau_m^+)\nonumber\\
&\qquad\quad=\frac{\partial
x(\tau_m^-,\btheta)}{\partial\theta_j}+\frac{\partial
\tau_m}{\partial\theta_j}f_x(\tau_m^-,\btheta)+\frac{\partial
\tau_m}{\partial \theta_j}\beta(\tau_m^+). \nonumber\end{align}
Since, at $\tau_m$ we have the event of discontinuity of
$\alpha_n$, and since by Assumption
\ref{assump:timeoutnd.nosimul}, $E_B$ cannot coincide with this
event, we have $\beta(\tau_m^-)=\beta(\tau_m^+)$. Using this in
the above expression by knowing that
$f_x(\tau_m^-,\btheta)=\alpha(\tau_m^-,\btheta)-\beta(\tau_m^-)$
reveals $$\frac{\partial h_1(m,\tau_k^-,\btheta)}{\partial
\theta_j}=\frac{\partial
x(\tau_m^-,\btheta)}{\partial\theta_j}+\tau_{m,j}^\prime\alpha(\tau_m^-,\btheta).
$$ Finally, notice that $h_1(m,\tau_k,\btheta)$ is not directly
dependent on the value of any state variable at $\tau_k$, hence
the first term in the paranthesized expression in
(\ref{induprime}) is 0. Putting all the results into
(\ref{induprime} proves part $(iv)$. Finally, for case $(v)$, by
(\ref{timer2}), we have
$h_2(m,t,\btheta)=x(\tau_m,\btheta)-(t-\tau_m)$ for
$t\in[\tau_m,\tau_k)$ with the boundary condition
$h_2(m,\tau_k,\btheta)=0$. Applying the same procedure proves part
$(v)$ and the whole lemma. \myqedhere

\subsubsection{State Derivatives}
Here, we find the derivative equations for the state variables $x$
and $\alpha$. We do not include other states in the IPA
calculations as by (\ref{delayedaprime}) and Lemma
\ref{lem:wprime}, they can be calculated in terms of the
derivatives of $x$ and $\alpha$:

\noindent\textbf{A) Analysis at event times}: We start by finding
the derivative update equation of $\alpha_n(t,\btheta)$. Notice
that unlike \cite{AliCassCDC11} where $\alpha_n(t,\btheta)$ was a
continuous function of time, here, by definition
$\alpha_n(t,\theta)$ has a discontinuity whenever event
$[w_n>\theta_n]$ occurs. In general, we have the following lemma:
%+[f_{\alpha_n}(\tau_k^-,\btheta)-f_{\alpha_n}(\tau_k^+,\btheta)]\tau_k^\prime
\begin{lemma} Concerning state variable $\alpha_n(t,\btheta)$,
\begin{equation}\frac{\partial\alpha_n(\tau_k^+,\btheta)}{\partial\theta_j}=0,
\agar e_k=[w>\theta_n].
\label{alphaprimediscontinuous}\end{equation} Also, for the rest
of the events, we have
\begin{align}&\frac{\partial\alpha_n(\tau_k^+,\btheta)}{\partial\theta_j}=\frac{\partial\alpha_n(\tau_k^+,\btheta)}{\partial\theta_j}\nonumber\\&\qquad+
\left\{\begin{array}{ll} -r_n\frac{\ind_{n=j}\tilde
\alpha(\tau_k,\btheta)-\frac{\partial
\tilde{x}(\tau_k^-,\btheta)}{\partial\theta_j}}{\tilde
\alpha(\tau_k^-,\btheta)-{\beta(\tau_k^-)}}
& \agar e_k=[w\leq \theta_n]\\
0 & \ow \end{array}\right. \label{alphaprimecontinuous}\end{align}
\label{lem:alphadiscupdate}
\end{lemma} \vspace{-3mm} \emph{Proof:} Equation \eqref{alphaprimediscontinuous} is immediate from
(\ref{ipaxdiscontdiscupdate}) noting $r(t)=\alpha_{n,\min}$ is
independent of $\btheta$.  Focusing on the second part, the proof
of this lemma is also straightforward by using
(\ref{ipaxcontdiscupdate}). Notice that the only event $e_k$
(excluding $[w>\theta_n]$) at whose occurrence time the difference
$f_{\alpha_n}(\tau_k^-,\btheta)-f_{\alpha_n}(\tau_k^+,\btheta)$ is
non-zero is $e_k=[w\leq \theta_n]$ where by (\ref{alphadot}), we
get
$f_{\alpha_n}(\tau_k^-,\btheta)-f_{\alpha_n}(\tau_k^+,\btheta)=0-r_n=-r_n$.
Using this information in (\ref{ipaxcontdiscupdate}) proves the
lemma. \myqedhere\\
Next, for any real-valued function $f$, we define
$\Delta\alpha_n(\tau_k)=\alpha_n(\tau_k^-,\btheta)-\alpha_n(\tau_k^+,\btheta)$
for $k=1,\ldots,K$. We derive the discrete update equations for
state variable $x$.
\begin{lemma} Concerning state variable $x(t,\btheta)$, we have
\begin{align}&\frac{\partial x(\tau_k^+,\btheta)}{\partial \theta_j}=\frac{\partial x(\tau_k^-,\btheta)}{\partial \theta_j}\nonumber\\&+\left\{\begin{array}{ll}
-\frac{\partial x(\tau_k^-,\btheta)}{\partial \theta_j} & \agar e_k=[x=0]\\
\Delta\alpha_n(\tau_k,\btheta)\frac{\ind_{n=j}\tilde
\alpha(\tau_k,\btheta)-\frac{\partial
\tilde{x}(\tau_k^-,\btheta)}{\partial\theta_j}}{\tilde\alpha(\tau_k^-,\btheta)-\beta(\tau_k^-)}
& \agar
e_k=[w>\theta_n]\\
0 & \ow
\end{array}\right.\label{xdiscupdate}\end{align}
\label{lem:xdiscupdate}
\end{lemma} \vspace{-3mm}
\emph{Proof:} The proof follows by straightforward application of
(\ref{ipaxcontdiscupdate}) to state variable $x$. When
$e_k=[x=0]$, we get
$f_x(\tau_k^-,\btheta)=\alpha(\tau_k^-,\btheta)-\beta(\tau_k^-)$
and $f_x(\tau_k^+,\btheta)=0$. Hence, by using $\frac{\partial
\tau_k}{\partial \theta_j}$ from Lemma \ref{lem:evtimederCh5}, we
find that
\begin{align}&[f_x(\tau_k^-,\btheta)-f_x(\tau_k^+,\btheta)]\frac{\partial
\tau_k}{\partial \theta_j}=\nonumber\\
&\quad
[\alpha(\tau_k^-,\btheta)-\beta(\tau_k^-)]\frac{-\frac{\partial
x(\tau_k^-,\btheta)}{\partial\theta_j}}{\alpha(\tau_k^-,\btheta)-\beta(\tau_k^-)}=-\frac{\partial
x(\tau_k^-,\btheta)}{\partial\theta_j},\nonumber\end{align} which
proves the first condition. When $e_k=[w>\theta_n]$, we again
invoke Lemma \ref{lem:evtimederCh5} for this event and use it in
(\ref{ipaxcontdiscupdate}) considering
$f_x(\tau_k^-,\btheta)-f_x(\tau_k^+,\btheta)=\alpha_n(\tau_k^-,\btheta)-\alpha(\tau_k^+,\btheta)$.
For other events, case by case analysis shows that
$f_x(\tau_k^-,\btheta)-f_x(\tau_k^+,\btheta)$ is only non-zero
when $E_B$ occurs. However, this is an exogenous event with
$\frac{\partial\tau_k}{\partial\theta_j}=0$ for all
$j=1,\ldots,N.$ \myqedhere

\noindent\textbf{B) Analysis between event times}: We start from
$\alpha_n(t,\btheta)$. Recall that by (\ref{alphadot}) and
Assumption \ref{assump:timeoutnd.infsupply},
$f_{\alpha_n}(t,\theta)=f_{\phi_n}(t,\theta)$ for all $n$ and
$t\in[0,T)$. Since $r_n$ is not a function of $\btheta$,
differentiating $f_{\alpha_n}(t,\theta)$ with respect to
$\theta_j$ reveals \begin{equation} \label{falphaprime}
\frac{\partial f_{\alpha_n}(t,\btheta)}{\partial \theta_j}=0,\quad
\forall\ t\in[0,T) \mbox{ and } j=1,\ldots,N.\end{equation}

For $x(t,\theta)$, (\ref{xndot}) implies that for all $t\in[0,T)$, $$f_x(t,\btheta)=\left\{\begin{array}{ll} 0 & \agar x(t,\theta)=0,\alpha(t,\btheta)\leq\beta(t)\\
                                                                       \alpha(t,\theta)-\beta(t)
                                                                       &
                                                                       \ow
                                                                      \end{array}\right. $$
therefore, for $j=1,\ldots,N$, \begin{equation}\label{fxprime} \frac{\partial f_x(t,\btheta)}{\partial\theta_j}=\left\{\begin{array}{ll} 0 & \agar x(t,\theta)=0,\alpha(t,\btheta)\leq\beta(t)\\
                                                                      \frac{\partial\alpha(t,\theta)}{\partial\theta_j} & \ow
                                                                    \end{array}\right.\end{equation}

Using (\ref{delayedaprime}) and Lemma \ref{lem:wprime} and
applying the IPA formulas given above completes the state
derivative estimations.

\subsubsection{IPA Implementation}  We assume that the
IPA derivatives for the states and event times over the interval
$[t-w(t,\btheta),t)$ are available for IPA calculations at $t$.

Recall that we are interested in estimating the derivative of the
average cost (\ref{J})  by finding the derivative of the sample
function (\ref{goodput}) which can be calculated from
(\ref{Gnprime}). However, evaluating (\ref{Gnprime}) is contingent
upon the knowledge of $\alpha_{n,j}^\prime(t,\btheta)$,
$\frac{d\alpha_n(t-\theta_n,\btheta)}{d\theta_j}$ and
$\tau_{k,j}^\prime$. By (\ref{delayedaprime}),
$\frac{d\alpha_n(t-\theta_n,\btheta)}{d\theta_j}$ can be readily
obtained given that the past values of
$\alpha^\prime_{n,j}(t,\btheta)$ and $w(t-\theta_n,\btheta)$ are
known. By Lemma \ref{lem:evtimederCh5}, $\tau_{k,j}^\prime$ can be
obtained if the derivatives
$\frac{\partial\alpha(t,\btheta)}{\partial\theta_j}$,
$\frac{\partial x(t,\btheta)}{\partial\theta_j}$ and
$\frac{\partial \tilde x(t,\btheta)}{\partial\theta_j}$ are
available. The latter is the value of $\frac{\partial
x(t,\btheta)}{\partial\theta_j}$ at $t-w(t,\btheta)$. The
derivative $\frac{\partial x(t,\btheta)}{\partial\theta_j}$ itself
can be found by evaluating
 (\ref{xdiscupdate}) at event times $[x=0]$ and $[w>\theta_n]$, $n=1,\ldots,N$ and integrating
 \eqref{fxprime} between these events. Furthermore, $\frac{\partial\alpha(t,\btheta)}{\partial\theta_j}$
 can be calculated by evaluating
\eqref{alphaprimediscontinuous} at the occurrences of
$[w>\theta_n]$ and \eqref{alphaprimecontinuous} at those of
$[w\leq\theta_n]$. Notice that by (\ref{falphaprime}), no
integration of $\frac{\partial
f_{\alpha_n}(t,\btheta)}{\partial\theta_j}$ is needed.

%\section{Numerical Results}\label{sec:timeoutnd_numericals}

\section{Conclusions}\label{sec:timeoutnd.conclusion}
We extended the results in \cite{AliCassCDC11} to multiple
communication links and laid down the conditions under which the
network buffer works according to the FCFS policy. We also showed
how extensions to non-FCFS policies are possible. Under the FCFS
policy, we derived the derivative estimates of average goodput at
each transmitter and used it to find the sensitivity of the
aggregate goodput with respect to timeout thresholds.

\bibliographystyle{abbrv}
\bibliography{AliBIB_FULL}

\end{document}